# Max/Min Puzzles in Geometry III†

Revised

James M. Parks
Dept. of Math.
SUNY Potsdam
Potsdam, NY
*parksjm@potsdam.edu*
February, 22, 2022

**Abstract**. The first two installments of this series of papers dealt with the maximum area polygons: Parallelograms, Rectangles, Squares or Equilateral Triangles, in given triangles. Minimum area polygons were also considered in the second paper on Equilateral Triangles. In this paper the puzzle will be turned the other way around. Given the regular unit polygons, a Square, Pentagon, or Hexagon, we searched for the smallest area triangle(s) which contains it.

The Dynamic Geometry Software *Sketchpad v5.10 BETA* was used for all figures in this paper.

**I**. The first regular polygon we will consider is the *unit square*. It will be assumed that triangles which enclose the square *inscribe it,* that is *the vertices of the square are all on the sides of the triangle*. It can be assumed that one side of the square is on a side of the triangle [4].

The following observation will be of use throughout the study of squares.
Let $\Delta ABC$ be a given triangle with base $AB$, in which a unit square $DEFG$ is inscribed, with $DE$ on side $AB,$ and assume there are no obtuse angles at vertices $A$ or $B,$ Fig. 1.
Then $\Delta ABC \sim \Delta GFC$, since the corresponding angles are equal [1].
Let $h$ be the height of $\Delta ABC$ at $C,$ $c = AB,$ $h'$ is the height of $\Delta GFC$ at $C$, and $c' = GF.$ Then $h/c = h'/c',$ but $h = (h' + 1),$ and $c' = 1,$ so $c = (h' + 1)/h'.$
Thus the following Lemma is established.

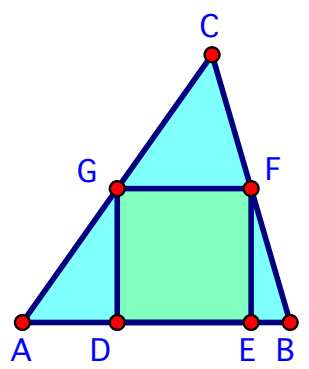

*Figure 1*

**Lemma 1**. *If $\Delta ABC$ inscribes the unit square DEFG with DE on AB, and AB has no obtuse angles at vertices A or B, then $AB = (h' + 1)/h'$, where $h'$ is the height of $\Delta GFC$ at C.*

**Examples 1**.

A. Let $\Delta ABC$ be the triangle which inscribes *DEFG,* as in Fig. 1, with $h' = 1$. Then, $c = 2$, and $(ABC) = 2\ u^2,$ by Lemma 1, where *(ABC) denotes the area of* $\Delta ABC$.

B. Let $\Delta ABC$ be the right angle triangle which inscribes *DEFG,* with sides of length *2*, and right angle at vertex *A* coinciding with the right angle in the unit square at vertex *D,* Fig. 2. Then $(ABC) = 2\ u^2$.

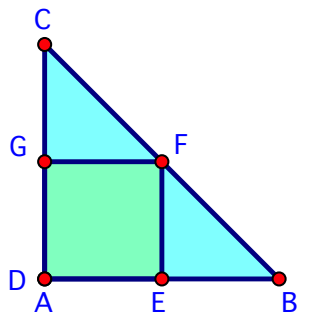

*Figure 2*

† Versions of some of these results have appeared in several settings previously [2], [3].

C. Let $\Delta ABC$ be an equilateral triangle which inscribes the unit square *DEFG,* Fig. 3. Then $\Delta GFC$ is also an equilateral triangle, and has base $GF = 1$. Thus the height of $\Delta GFC$ is $h' = \sqrt{3}/2$, and the height of $\Delta ABC$ is $h = (1 + \sqrt{3}/2)$. Hence, $c = (2\sqrt{3} + 3)/3$, by Lemma 1, so $(ABC) = (12+7\sqrt{3})/12\ u^2 \sim 2.01\ u^2$.

D. Let $\Delta ABC$ be the isosceles triangle which inscribes the unit square *DEFG,* with $CA = CB$, and central angle $\angle C = 120°$, Fig. 4. Then $\angle A = 30° = \angle B$, and $AD = EB = \sqrt{3}$. So $c = AB = (1 + 2\sqrt{3})$, and $h = (\sqrt{3}/6 + 1)$, Lemma 1. Thus, $(ABC) = (1 + 2\sqrt{3})(\sqrt{3}/6 + 1)/2\ u^2 \sim 2.88\ u^2$.

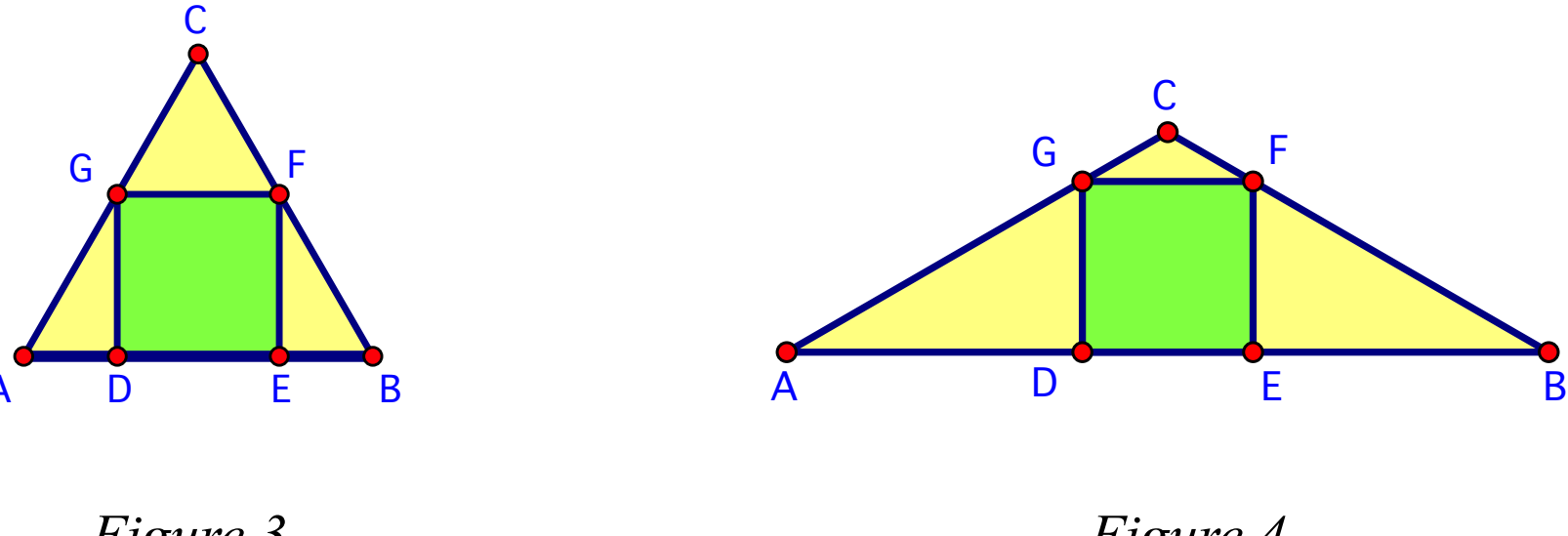

*Figure 3* *Figure 4*

The above examples lead us to conjecture that $2\ u^2$ is a lower bound for the minimum area of a triangle which inscribes a unit square. We state the conjecture as Theorem 1, a proof follows.

**Theorem 1**. *If a unit square DEFG is inscribed in* $\Delta ABC$, *with DE on side AB, then the height h of* $\Delta ABC$ *from C satisfies* $h = 2$, *iff* $\min(ABC) = 2\ u^2$.

By Lemma 1, if $\Delta ABC$ is a triangle in which a unit square *DEFG* is inscribed, and $\Delta ABC$ has height $h = 2$, measured from *C,* then the base $c = 2$, so $(ABC) = 2\ u^2$, Fig. 1.
Suppose $\Delta RST$ contains the inscribed unit square *DEFG,* with *DE* on side *RS,* and has height $h \neq 2$. For instance $h = (2 + e)$, for $e > 0$. By Lemma 1, $h' = (1 + e)$, so $c = (2 + e)/(1 + e)$, and $(RST) = (2 + e)^2/(2 + 2e)$. It then follows that $(RST) > 2\ u^2$, since the inequality $(2 + e)^2/(2 + 2e) > 2\ u^2$ is equivalent to $e^2 > 0$, which is obviously true.
A similar argument holds if $\Delta RST$ has height $h = (2 - e)$, for $0 < e < 1$.
Thus the Theorem is established.

**II**. The next polygon to be considered is the *regular unit pentagon,* Fig.5.
The interior angles of pentagon *DEFGH* are $108°$, since a pentagon interior can be divided into *3* triangles.

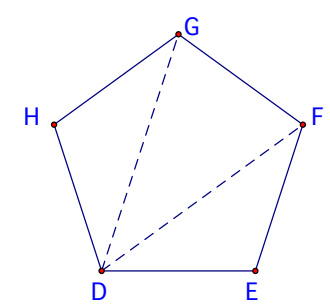

*Figure 5*

The major chords, such as *HF*, all have length $(1 + \sqrt{5})/2$, the Golden Ratio. This follows, since the chord *HF,* with length *s,* determines an isosceles triangle, *FGH,* with base angles $36°$. It also determines the isosceles trapezoid *DEFH,* with sides of length *1, 1, 1,* and *s.* Since $<HFG = <DFH = <EFD = 36°$, and the diagonals all have length s, the interior triangles *DFH* and *EFH* are isosceles with base angles of $72°$, Fig. 6. So the equation $s^2 = s + 1$ results, by Ptolemy's Theorem. The positive solution is $s = (1 + \sqrt{5})/2$.

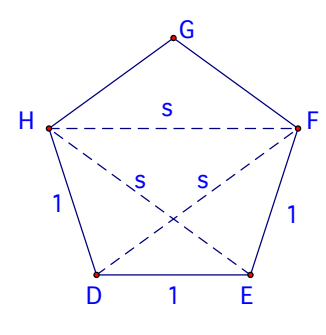

*Figure 6*

Another geometric argument is given in Theorem 2.

It is only possible to inscribe a pentagon in a triangle in a few ways. So it will be necessary to also consider triangles which 'enclose' a pentagon in a more general way, see Fig 7. This requires that 4 of the 5 vertices of the pentagon must be on the enclosing triangle.

***It will be assumed that there are no obtuse angles on the base side AB of ΔABC from here on, unless noted otherwise***.

If the line segment *HF* is added to the unit regular pentagon *DEFGH*, then any triangle Δ*ABC,* which encloses the pentagon such that *DE* is on side *AB, F* is on *CB,* and *H* is on *CA,* is similar to Δ*HFC,* by angle equality, Fig.7.
Thus, if $h$ is the height of Δ*ABC* from *C, h'* is the height of Δ*HFC* from *C*, $s = HF$, $c = AB$, and $d$ is the height of *DEFH,* then $h/c = h'/s$. But $h = h' + d$, so $c = (h' + d)s/h'$.
This proves Lemma 2.

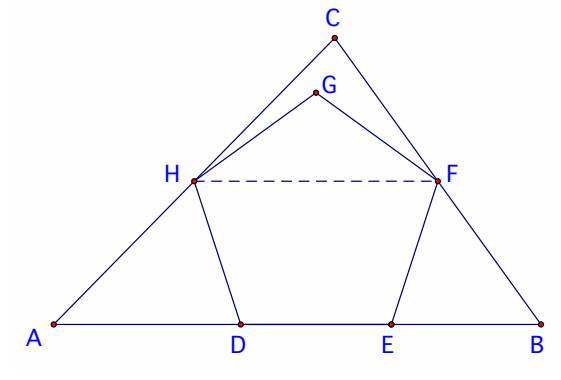

*Figure 7*

**Lemma 2.** *Let* Δ*ABC enclose the unit regular pentagon DEFGH, with DE on AB, F on CB, and H on CA. Also, assume h is the height of* Δ*ABC from C, h' is the height of* Δ*HFC from C,* $c = AB$, $s = HF$, *and d* is the height of *DEFH. Then* $c = (h' + d)s/h'$.

**Examples 2**.

A. In Fig. 8, the triangle Δ*ABC* which encloses the regular unit pentagon is special. Using *Sketchpad* the height is $h = \sqrt{((5 + \sqrt{5})/2)} \sim 1.90$, which is twice the height of the line segment *HF* above the base *AB* (see Th. 2 below for details on this). Then the base $AB = 2s = (1 + \sqrt{5}) \sim 3.24$, where $s = HF$, and $(ABC) = \sqrt{((5 + \sqrt{5})/2)}(1 + \sqrt{5})/2 \sim 3.08\ u^2$.

B. Let Δ*ABC* be the triangle which encloses the regular unit pentagon *DEFGH,* with side *DE* on *AB,* and height $h = 2$, Fig. 9. Then, by Lemma 2 and *Sketchpad*, the base $c = AB \sim 3.09$, and $(ABC) \sim 3.08\ u^2$.

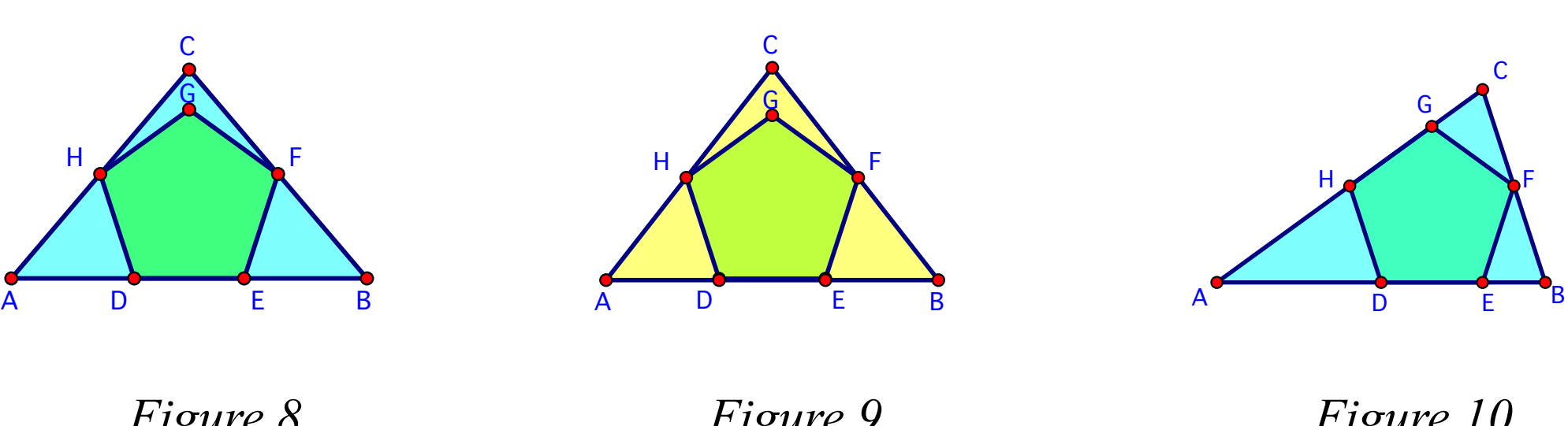

*Figure 8* *Figure 9* *Figure 10*

C. Let Δ*ABC* be the triangle in which the regular unit pentagon *DEFGH* is inscribed*,* with *DE* on side *AB,* where *G* is on side *CA,* and the sides *CH* and *ED* are extended to intersect at *A,* and the side $BC$ is perpendicular to the line segment $AF$, Fig. 10. It then follows, that $(ABC) \sim 3.08\ u^2$, by *Sketchpad* (see also Theorem 2).

D. Let $\Delta ABC$ be the triangle in which the regular unit pentagon $DEFGH$ is inscribed, with $DE$ on side $AB$, the vertex $C$ coincides with vertex $G$, and the sides $CH$ and $CF$ are extended to intersect with the line on the base side at $A$ and $B$, respectively.
Then $\angle A = \angle B = 36°$, and $\angle C = 108°$.
Using *Sketchpad* it follows that $(ABC) \sim 3.26\ u^2$.

*Figure 11*

These examples lead us to conjecture Theorem 2.

**Theorem 2**. *Let DEFGH be a regular unit pentagon, and let the side DE be reflected about chord line HF to form the line segment ST. Let $\Delta ABC$ enclose the pentagon subject to the following conditions: side DE is on side AB of the triangle, there are no obtuse angles on AB, and assume vertex F is on BC, vertex H is on CA, and vertex C is on ST. Then $\Delta ABC$ is a minimum enclosing triangle of the pentagon with $(ABC) = \sqrt{5+2\sqrt{5}}\ u^2 \sim 3.08\ u^2$.*

Proof. Construct the triangle $\Delta ABC$ as in Example $C$, Fig.10. Since the exterior angles of the pentagon are $72°$, the triangles $\Delta ABC$, $\Delta ADH$, $\Delta BFE$ and $\Delta CGF$ are all similar isosceles triangles with base angles of $72°$, and central angles of $36°$.
Also, a copy of $\Delta ABC$, labeled $\Delta A'B'C'$, can be constructed to the right of the pentagon with the same properties as $\Delta ABC$, see dashed lines in Fig. 12.

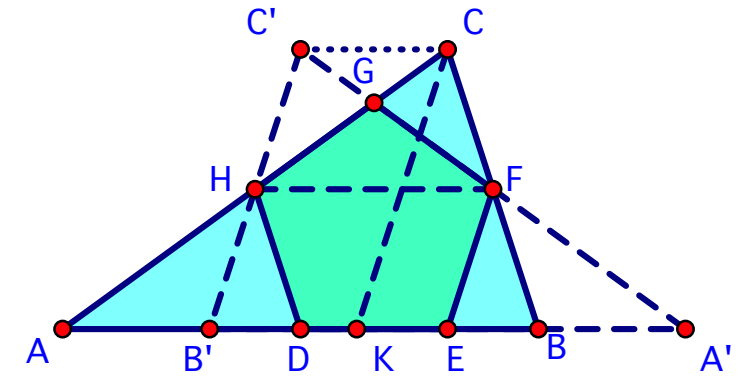

*Figure 12*

Since $H$ and $F$ are the midpoints of $B'C'$ and $BC$, chord $HF$ is halfway between $DE$ and $C'C$, and $BC = B'C' = 2$. This also means that $C'C$ coincides with $ST$, the reflection of $DE$ about $HF$.
Since $\Delta HFC \sim \Delta ABC$, it follows that $2/1 = (s+1+e)/s$, where $s = HF = HC = AD$, $e = EB$, and $c = AB = s + 1 + e$. So $s = 1 + e$, and $c = AB = 2s$.
Let $K$ be the point on $AB$ such that $\Delta KBC$ is an isosceles triangle, Fig. 12.
Then $KB = 2e$, so $2s/2 = 2/2e$, or $s = 1/e$, and $s = 1/(s-1)$.
Thus, since $s = (1 + \sqrt{5})/2$, we have $e = 1/s = (\sqrt{5} - 1)/2$, and $AB = 1 + \sqrt{5}$.
Hence, $\Delta EBC$ is a right triangle, since $\angle B = 72°$, and $\angle ECB = 18°$, from the observation above that $KB = 2e$.
So, if $h = EC$, then $h^2 + e^2 = 4$, and $h = \sqrt{(5 + \sqrt{5})/2} \sim 1.90$.
It then follows that $(ABC) = hc/2 = \sqrt{5 + 2\sqrt{5}}\ u^2 \sim 3.08\ u^2$.
If $\Delta ABC$ is a 3rd triangle which encloses the pentagon $DEFGH$ and has vertex $C''$ on $ST = C'C$, then $\Delta HFC \sim \Delta ABC$, and the ratios satisfy $h/h' = c/c'$, where $c = AB$, and $c' = HF$, Fig. 13.
But $h = 2h'$, and $c = 2c'$, where $c' = s$, so $(ABC) = hc/2 = \sqrt{5 + 2\sqrt{5}}\ u^2 \sim 3.08\ u^2$.

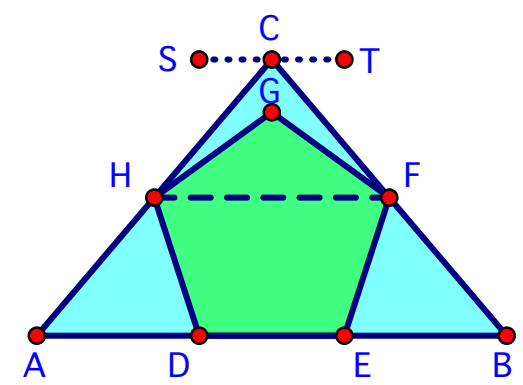

*Figure 13*

Now suppose $\Delta A'B'C'$ is another triangle which encloses pentagon *DEFGH,* with *DE* on *A'B',* and there are no obtuse angles on the base *A'B',* where vertex *F* is on *B'C',* vertex *H* is on *C'A',* and vertex *C'* is <u>not</u> on segment *ST,* Fig. 14. Then, for $c' = A'B'$, and $h'$ the height of *C'* over *A'B',* $h'/c' = (h' - h/2)/s$, by Lemma 2, where *h* is the distance from *ST* to *A'B',* $s = HF$.
Thus $(A'B'C) = sh'^2/(2h' - h)$.

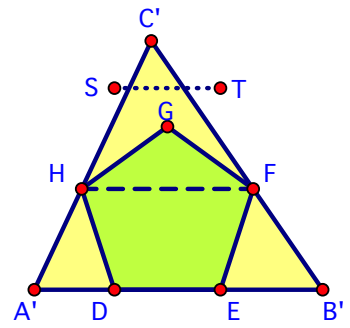

*Figure 14*

The area of the first triangle above in Fig. 13, with vertex *C* on *ST* is *sh,* since $AB = 2s$. So the claim is that the area of the triangle in Fig. 14 is greater than the area of the triangle in Fig. 13, that is $sh'^2/(2h' - h) > sh$.
But this is equivalent to the inequality $(h' - h)^2 > 0$, which is obviously true, regardless of whether $h'$ is greater than, or less than *h*.
Hence the claim is true.

**III**. The last polygon we will considered is the *regular unit hexagon,* Fig.15. This polygon is special, since, as we will see, it does not follow the patterns of the areas for the square and the pentagon above.

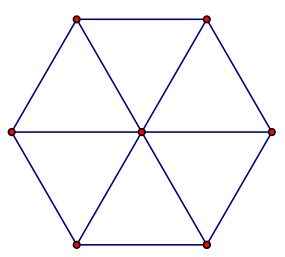
*Figure 15*

The interior angles of the hexagon are $120°$, hence the exterior angles are $60°$, since the regular hexagon can be divided into *6* unit equilateral triangles.
So the height of the unit hexagon is $\sqrt{3}$, and since the equilateral triangles each have area $\sqrt{3}/4\ u^2$, the area of the unit hexagon is $3\sqrt{3}/2\ u^2$.

**Examples 3**.

A. If $\Delta ABC$ is the equilateral triangle which inscribes the unit regular hexagon *DEFGHJ,* Fig. 16. Then, the height *h* of $\Delta ABC$ is $3\sqrt{3}/2$, and it's sides have length *3,* so $(ABC) = 9\sqrt{3}/4\ u^2 \sim 3.90\ u^2$.

B. Let $\Delta ABC$ be the yellow isosceles triangle which encloses the unit regular hexagon *DEFGHJ*, and has height *3,* and base $AB \sim 2.81$, Fig. 17. By *Sketchpad* the approximate area of the triangle is $(ABC) \sim 4.22\ u^2$.

C. Let $\Delta A'B'C'$ be the blue triangle which encloses the unit regular hexagon *DEFGHJ,* and has height *2.25,* and base $A'B' \sim 4.34$, Fig. 17. By *Sketchpad* the approximate area of the triangle is $(A'B'C') \sim 4.89\ u^2$.

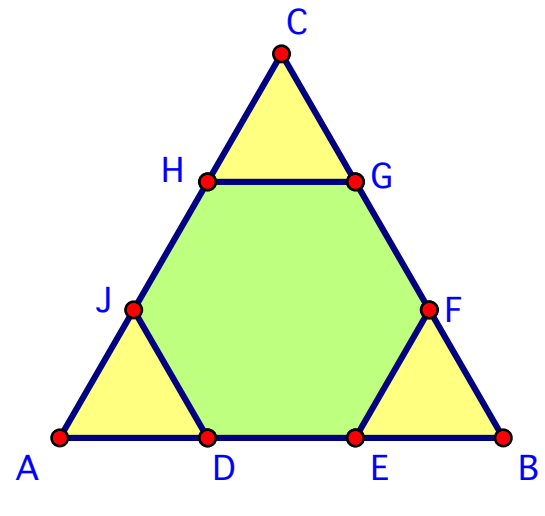

*Figure 16*

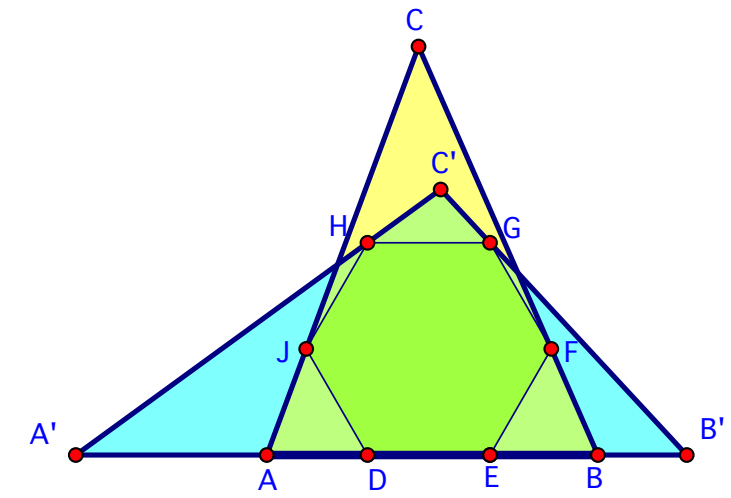

*Figure 17*

Using these Examples as a guide we can conjecture Theorem 3.

**Theorem 3**. *The minimum triangle* $\Delta ABC$ *in area which encloses the unit regular hexagon is an equilateral triangle with sides of length 3, and* $(ABC) = 9\sqrt{3}/4\ u^2 \sim 3.90\ u^2$.

Proof. If $\Delta ABC$ is the equilateral triangle in which is encloses the unit regular hexagon *DEFGHJ,* then $(ABC) = 9\sqrt{3}/4\ u^2$, Example 3A, Fig. 16.
Let $\Delta A'B'C'$ be another triangle which encloses the unit regular hexagon, with height $h' > 3\sqrt{3}/2 \sim 2.60$, Fig. 18.
We claim that $(A'B'C') > (ABC)$.
Connect vertices $F$ and $J$ for the chord $FJ$, then $\Delta JFC' \sim \Delta A'B'C'$, by equal angles.
Let $c' = A'B'$, and note $JF = 2$, then $h'/c' = (h' - \sqrt{3}/2)/2$.
Hence $c' = 2h'/(h' - \sqrt{3}/2)$, and $(A'B'C') = h'^2/(h' - \sqrt{3}/2)\ u^2$.
Thus $h'^2/(h' - \sqrt{3}/2) > 9\sqrt{3}/4$ iff $h'^2 - 9\sqrt{3}h'/4 + 27/8 > 0$.
But $h'^2 - 9\sqrt{3}h'/4 + 27/8 = (h' - 3\sqrt{3}/4)(h' - 3\sqrt{3}/2)$, and this product is positive whenever $h' > 3\sqrt{3}/2$.
Hence the inequality holds.
A similar argument works for $\sqrt{3} < h' < 3\sqrt{3}/2$.
Thus the conjecture is true.

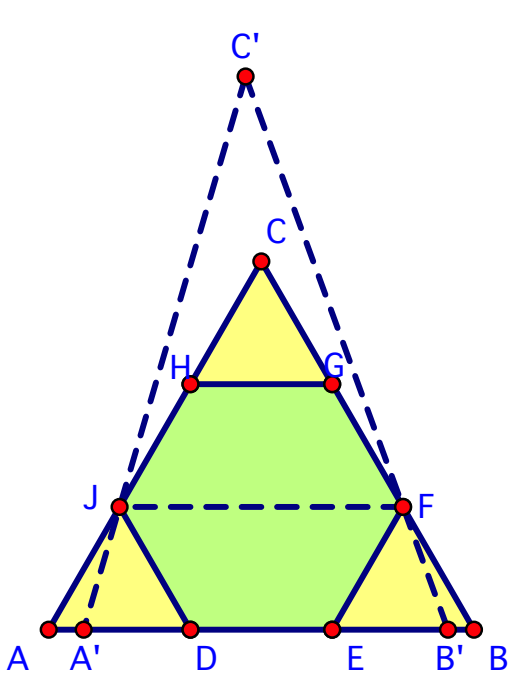

*Figure 18*

The general case for regular unit $n$-gons, the regular unit polygons with $n$ sides, is considered in [3] for $n \geq 4$. In that paper the patterns of the minimum enclosing triangles which appear above are shown to hold in the more general case for all regular unit $n$-gons, $n \geq 4$.
If we use the notation $n = r\bmod 3$, and assume $r = 0$, then it is demonstrated in [3] that these $n$-gons have minimum enclosing triangles which are equilateral triangles, as in the case of the regular unit hexagon, a $6$-gon, in Section III above.
If $r \neq 0$, then the unit $n$-gons have minimum enclosing triangles similar to the cases for the unit square, a $4$-gon, Section I, and the regular unit pentagon, a $5$-gon, Section II.

**References.**